\documentclass[12pt]{article}
\usepackage{amssymb,amsmath,amsfonts,t1enc}
\newcommand{\Z}{\mathbb Z}
\newcommand{\Q}{\mathbb Q}
\newcommand{\R}{\mathbb R}
\newcommand{\C}{\mathbb C}
\newcommand{\D}{\textbf{\textit{D}}}
\newcommand{\G}{\textbf{\textit{G}}}
\newcommand{\K}{\textbf{\textit{K}}}
\newcommand{\M}{\textbf{\textit{M}}}
\newcommand{\LL}{\textbf{\textit{L}}}

\begin{document}
\thispagestyle{empty}
\hfill
\par
\noindent
\centerline{{\large Arithmetic neighbourhoods of numbers}}
\vskip 0.2truecm
\centerline{{\large Apoloniusz Tyszka}}
\vskip 0.2truecm
\par
\noindent
{\bf Abstract.} Let $\K$ be a ring and let $A$ be a subset of $\K$.
We say that a map $f:A \to \K$ is {\sl arithmetic}
if it satisfies the following conditions:
if $1 \in A$ then $f(1)=1$,
if $a,b \in A$ and $a+b \in A$ then $f(a+b)=f(a)+f(b)$,
if $a,b \in A$ and $a \cdot b \in A$ then $f(a \cdot b)=f(a) \cdot f(b)$.
We call an element $r \in \K$ {\sl arithmetically fixed} if there is
a finite set $A \subseteq \K$ (an {\sl arithmetic neighbourhood} of $r$
inside~$\K$) with $r \in A$ such that each arithmetic map $f:A \to \K$
fixes $r$, i.e. $f(r)=r$. We prove:
for infinitely many integers~$r$ for some arithmetic neighbourhood of~$r$
inside~$\Z$ this neighbourhood is a neighbourhood of~$r$ inside~$\R$
and is not a neighbourhood of~$r$ inside~${\Z}[\sqrt{-1}]$;
for infinitely many integers~$r$ for some arithmetic neighbourhood of~$r$
inside~$\Z$ this neighbourhood is not a neighbourhood of~$r$ inside~$\Q$;
if $\K={\Q}(\sqrt{5})$ or $\K={\Q}(\sqrt{33})$, then
for infinitely many rational numbers~$r$ for some
arithmetic neighbourhood of~$r$ inside~$\Q$ this neighbourhood
is not a neighbourhood of~$r$ inside~$\K$;
for each $n \in (\Z \cap [3,\infty)) \setminus \{2^2,2^3,2^4,...\}$
there exists a finite set~${\cal J}(n) \subseteq \Q$
such that ${\cal J}(n)$ is a neighbourhood of~$n$ inside~$\R$ and ${\cal J}(n)$
is not a neighbourhood of~$n$ inside~$\C$.
\vskip 0.2truecm
\par
Let $\K$ be a ring and let $A$ be a subset of $\K$.
We say that a map $f:A \to \K$ is {\sl arithmetic}
if it satisfies the following conditions:
\par
\noindent
{\bf (1)}~~~if $1 \in A$ then $f(1)=1$,
\par
\noindent
{\bf (2)}~~~if $a,b \in A$ and $a+b \in A$ then $f(a+b)=f(a)+f(b)$,
\par
\noindent
{\bf (3)}~~~if $a,b \in A$ and $a \cdot b \in A$ then $f(a \cdot b)=f(a) \cdot f(b)$.
\par
\noindent
We call an element $r \in \K$ {\sl arithmetically fixed} if there is
a finite set $A \subseteq \K$ (an {\sl arithmetic neighbourhood} of $r$ inside~$\K$)
with $r \in A$ such that each arithmetic map $f:A \to \K$ fixes $r$,
i.e. $f(r)=r$.
All previous articles on arithmetic neighbourhoods (\cite{Tyszka1}, \cite{Lettl}, \cite{Tyszka2})
dealt with a description of a situation where for an element in a field there exists
an arithmetic neighbourhood. If $\K$ is a field, then any $r \in \K$ is arithmetically
fixed if and only if $\{r\}$ is existentially first-order definable in the language
of rings without parameters (\cite{Tyszka2}). Therefore, presentation of the arithmetic
neighbourhood of the element $r$ belonging to the field $\K$ is the simplest way of
expression that $\{r\}$ is existentially definable in $\K$.
\vskip 0.2truecm
\par
We want to find integers~$r$ with property~{\bf (4)},
integers~$r$ with property~{\bf (5)}, and
rational numbers~$r$ with property~{\bf (6)}.
\par
\noindent
{\bf (4)}~~Each arithmetic neighbourhood of~$r$ inside~$\Z$ is also a neighbourhood of~$r$ inside each ring extending~$\Z$.
\par
\noindent
{\bf (5)}~~Each arithmetic neighbourhood of~$r$ inside~$\Z$ is also a neighbourhood of~$r$ inside~$\Q$.
\par
\noindent
{\bf (6)}~~Each arithmetic neighbourhood of~$r$ inside~$\Q$ is also a neighbourhood of~$r$ inside each ring
extending~$\Q$.
\vskip 0.2truecm
\par
\noindent
Obviously, condition~{\bf (4)} implies condition~{\bf (5)}.
\vskip 0.2truecm
\par
By condition~{\bf (1)} for any ring~$\K$
each arithmetic neighbourhood of~$1$ inside~$\K$
is also a neighbourhood of~$1$ inside each ring extending~$\K$.
Since~$0+0=0$, by condition~{\bf (2)} for any ring~$\K$
each arithmetic neighbourhood of~$0$ inside~$\K$
is also a neighbourhood of~$0$ inside each ring extending~$\K$.
\vskip 0.2truecm
\par
We prove that for any ring~$\K$ with~$2 \neq 0$
each arithmetic neighbourhood of~$2$ inside~$\K$
is also a neighbourhood of~$2$ inside each ring~$\LL$ extending~$\K$.
Assume that $A$ is an arithmetic neighbourhood of~$2$ inside~$\K$
and $f:A \to \LL$ is an arithmetic map.
Then $1 \in A$, because in the opposite case the arithmetic map $A \to \{0\}$ moves~$2$,
which is impossible. Since $f$ satisfies conditions~{\bf (1)} and~{\bf (2)}, we get
$f(2)=f(1+1)=f(1)+f(1)=1+1=2$.
\vskip 0.2truecm
We prove that for any ring~$\K$ with~$\frac{1}{2} \in \K$
each arithmetic neighbourhood of~$\frac{1}{2}$ inside~$\K$
is also a neighbourhood of~$\frac{1}{2}$ inside each ring~$\LL$ extending~$\K$.
Assume that $A$ is an arithmetic neighbourhood of~$\frac{1}{2}$ inside~$\K$
and $f:A \to \LL$ is an arithmetic map.
Then $1 \in A$, because in the opposite case
the arithmetic map $A \to \{0\}$ moves~$\frac{1}{2}$, which is impossible.
Since $f$ satisfies conditions~{\bf (1)} and~{\bf (2)}, we get
$1=f(1)=f(\frac{1}{2}+\frac{1}{2})=f(\frac{1}{2})+f(\frac{1}{2})$.
Hence, $f(\frac{1}{2})=\frac{1}{2}$.
\vskip 0.2truecm
\par
The above results imply that the numbers $r=1$, $r=0$, $r=2$
satisfy conditions~\mbox{{\bf (4)}--{\bf (6)},} and $r=\frac{1}{2}$
satisfies condition~{\bf (6)}.
\vskip 0.2truecm
\par
Let $\K$ be an algebraically closed field and $r \in \K$ is
arithmetically fixed. Then $r$ belongs
to the prime field in $\K$, see~\cite{Tyszka2}, cf.~\cite{Tyszka1}.
Let $A=\{x_1,...,x_n\}$ be a neighbourhood of~$r$,
$x_i \neq x_j$ if $ i \neq j$, and $x_1=r$.
We choose all formulae
$x_i=1$ ($1 \leq i \leq n$), $x_i+x_j=x_k$, $x_i \cdot x_j=x_k$
($1 \leq i \leq j \leq n$, $1 \leq k \leq n$) that are satisfied
in~$A$. Joining these formulae with conjunctions we get
some formula~$\Phi$.
Let $V$ denote the set of variables in $\Phi$.
Since $A$ is a neighbourhood of $r$ inside $\K$, we have
$$
\K \models \underbrace{...~\forall x_s~...}_{x_s \in \{x_1\} \cup V}~(\Phi \Rightarrow x_1=r)
$$
Of course, $\{x_1\} \cup V=V$ ([15, the proof of Theorem~2] and [16, the proof of Theorem~1])
but this equality will not be used later.
\vskip 0.2truecm
\par
\noindent
{\bf Proposition 1.} Let $\K$ be an algebraically closed field and
$r \in \K$ is arithmetically fixed. Then each arithmetic neighbourhood
of~$r$ inside~$\K$ is also a neighbourhood of~$r$ inside each integral
domain~$\D$ extending~$\K$.
\vskip 0.2truecm
\par
\noindent
{\it Proof.} We give a model-theoretic proof,
an alternative proof follows from Hilbert's Nullstellensatz.
Let $A$ be a neighbourhood of $r$ inside $\K$.
Let ${\D}_1$ denote the algebraic closure of the fraction field
of $\D$. It suffices to prove that $A$ is a neighbourhood
of~$r$ inside~${\D}_1$. Since $\K$ is a subfield of~${\D}_1$ and
every embedding between algebraically closed fields is
elementary ([6,~pp.~103~and~57], we obtain
$$
{\D}_1 \models \underbrace{...~\forall x_s~...}_{x_s \in \{x_1\} \cup V}~(\Phi \Rightarrow x_1=r)
$$
It implies that $A$ is a neighbourhood of~$r$ inside~${\D}_1$.
\newline
\rightline{$\Box$}
\vskip 0.2truecm
\par
Let ${\cal T}$ denote the elementary theory of integral domains
of characteristic $0$.
\par
\noindent
{\bf Proposition 2}. Let $\K$ be an algebraically closed field that extends $\Q$ and $r \in \K$
is arithmetically fixed.
\par
\noindent
{\bf (I)}~~~If $r=0$, then
$$
{\cal T} \vdash \underbrace{...~\forall x_s~...}_{x_s \in \{x_1\} \cup V}~(\Phi \Rightarrow x_1=0)
$$
{\bf (II)}~~~If $r=\frac{k}{w}$ for some $k,w \in \{1,2,3,...\}$, then
$$
{\cal T} \vdash \underbrace{...~\forall x_s~...}_{x_s \in \{x_1\} \cup V}~(\Phi \Rightarrow \underbrace{(1+...+1)}_{{\rm w-times}} \cdot x_1=\underbrace{1+...+1}_{{\rm k-times}})
$$
{\bf (III)}~~~If $r=-\frac{k}{w}$ for some $k,w \in \{1,2,3,...\}$, then
$$
{\cal T} \vdash \underbrace{...~\forall x_s~...}_{x_s \in \{x_1\} \cup V}~(\Phi \Rightarrow \underbrace{(1+...+1)}_{{\rm w-times}} \cdot x_1+\underbrace{1+...+1}_{{\rm k-times}}=0)
$$
{\it Proof.} We prove {\bf (II)} and omit similar proofs of {\bf (I)} and {\bf (III)}.
It suffices to prove that the sentence
$$
\underbrace{...~\forall x_s~...}_{x_s \in \{x_1\} \cup V}~(\Phi \Rightarrow \underbrace{(1+...+1)}_{{\rm w-times}} \cdot x_1=\underbrace{1+...+1}_{{\rm k-times}})
$$
holds true in each integral domain of characteristic $0$. Let $\G$ be any
integral domain of characteristic~$0$. Let ${\G}_1$ denote the 
algebraic closure of the fraction field of $\G$.
There exists an algebraically closed
field $\M$ such that both $\K$ and ${\G}_1$ embed into~$\M$. Of course,
$$
\K \models \underbrace{...~\forall x_s~...}_{x_s \in \{x_1\} \cup V}~(\Phi \Rightarrow \underbrace{(1+...+1)}_{{\rm w-times}} \cdot x_1=\underbrace{1+...+1}_{{\rm k-times}})
$$
Since every two algebraically closed fields of the same characteristic
are elementary equivalent [6,~p.~57], we obtain
$$
\M \models \underbrace{...~\forall x_s~...}_{x_s \in \{x_1\} \cup V}~(\Phi \Rightarrow \underbrace{(1+...+1)}_{{\rm w-times}} \cdot x_1=\underbrace{1+...+1}_{{\rm k-times}})
$$
Since $\G$ embeds into $\M$, we obtain
$$
\G \models \underbrace{...~\forall x_s~...}_{x_s \in \{x_1\} \cup V}~(\Phi \Rightarrow \underbrace{(1+...+1)}_{{\rm w-times}} \cdot x_1=\underbrace{1+...+1}_{{\rm k-times}})
$$
\newline
\rightline{$\Box$}
\vskip 0.2truecm
\par
Let $n \in \Z$, $n \geq 3$, $S_n=\{1,10,20,30\} \cup \{3,3^2,3^3,...,3^n\}$, $S=\bigcup_{n=3}^\infty\limits S_n$.
\vskip 0.2truecm
\par
\noindent
{\bf Theorem 1.} There is an arithmetic map $\gamma: S \to {\Z}[\sqrt{-1}]$
which moves all \mbox{$r \in S \setminus \{1\}$.} For each $r \in S_n \setminus \{1\}$
we have:
\par
\noindent
{\bf (7)}~~~$S_n$ is an arithmetic neighbourhood of~$r$ inside~$\R$, and so too inside~$\Q$ and~$\Z$,
\par
\noindent
{\bf (8)}~~~$S_n$ is not an arithmetic neighbourhood of~$r$ inside~${\Z}[\sqrt{-1}]$.
\vskip 0.2truecm
\par
\noindent
{\it Proof.} We prove {\bf (7)}. Assume that $f:S_n \to \R$ is an arithmetic map. Then,
\par
\noindent
$f(1)=1$,
\par
\noindent
$f(9)=f(3 \cdot 3)=f(3) \cdot f(3)=(f(3))^2$,
\par
\noindent
$f(27)=f(3 \cdot 9)=f(3) \cdot f(9)=f(3) \cdot (f(3))^2=(f(3))^3$,
\par
\noindent
$f(30)=f(27+3)=f(27)+f(3)=(f(3))^3+f(3)$,
\par
\noindent
$f(10)=f(9+1)=f(9)+f(1)=(f(3))^2+1$,
\par
\noindent
$f(20)=f(10+10)=f(10)+f(10)=(f(3))^2+1+(f(3))^2+1=2 \cdot (f(3))^2+2$,
\par
\noindent
$f(30)=f(20+10)=f(20)+f(10)=2 \cdot (f(3))^2+2+(f(3))^2+1=3 \cdot (f(3))^2+3$.
\par
\noindent
Therefore, $(f(3))^3+f(3)=3 \cdot (f(3))^2+3$. Hence $(f(3)-3) \cdot ((f(3))^2+1)=0$.
Thus $f(3)=3$, and by induction we obtain $f(3^k)=3^k$ for each $k \in \{1,2,3,...,n\}$.
Consequently,
\par
\noindent
$f(10)=f(9+1)=f(9)+f(1)=9+1=10$,
\par
\noindent
$f(20)=f(10+10)=f(10)+f(10)=10+10=20$,
\par
\noindent
$f(30)=f(20+10)=f(20)+f(10)=20+10=30$.
\par
\noindent
We have proved {\bf (7)}.
We define $\gamma:S \to {\Z}[\sqrt{-1}]$ as
\par
\noindent
\centerline{$\left\{(1,1),~(10,0),~(20,0),~(30,0)\right\} \cup
\left\{(3,\sqrt{-1}),~(3^2,(\sqrt{-1})^2),~(3^3,(\sqrt{-1})^3),~...\right\}$}
\par
\noindent
The map $\gamma$ is arithmetic and $\gamma$ moves all $r \in S \setminus \{1\}$,
so condition~{\bf (8)} holds true.
\newline
\rightline{$\Box$}
\par
We state a similar result without a proof. Let $n \in \Z$, $n \geq 1$,
\par
\noindent
$T_n=\{-2,1,5,10,20\} \cup \{4^1,...,4^n\}$, $T=\bigcup_{n=1}^\infty\limits T_n$.
Let us define $\tau: T \to {\Z}[\sqrt{-1}]$ as
\par
\noindent
\centerline{$\left\{(-2,\sqrt{-1}),~(1,1),~(5,0),~(10,0),~(20,0)\right\}
\cup \left\{(4,-1),~(4^2,1),~(4^3,-1),~(4^4,1),~...\right\}$}
\par
\noindent
The map $\tau$ is arithmetic and $\tau$ moves all $r \in T \setminus \{1\}$.
For each $r \in T_n \setminus \{-2,1\}$ we have:
\par
\noindent
$T_n$ is an arithmetic neighbourhood of~$r$ inside~$\R$, and so too inside~$\Q$ and~$\Z$,
\par
\noindent
$T_n$ is not an arithmetic neighbourhood of~$r$ inside~${\Z}[\sqrt{-1}]$.
\vskip 0.2truecm
\par
\noindent
{\bf Remark 1.} By Theorem~1 for infinitely many integers~$r$ fail both conditions~{\bf (4)} and~{\bf (6)}.
In Theorems~5, 6, and~7 we describe some other rational numbers~$r$ without property~{\bf (6)}.
\vskip 0.2truecm
\par
Let $n \in \Z$, $n \geq 3$, $B_n=\{1,5,25,26\} \cup \{3,3^2,3^3,...,3^n\}$,~$B=\bigcup_{n=3}^\infty\limits B_n$.
\vskip 0.2truecm
\par
\noindent
{\bf Theorem 2.} There is an arithmetic map $\phi: B \to \Q$
which moves all $r \in B \setminus \{1\}$.
For each $r \in B_n \setminus \{1,5\}$ we have:
\par
\noindent
{\bf (9)}~~~$B_n$ is an arithmetic neighbourhood of~$r$ inside~$\Z$,
\par
\noindent
{\bf (10)}~~~$B_n$ is not an arithmetic neighbourhood of~$r$ inside~$\Q$.
\vskip 0.2truecm
\par
\noindent
{\it Proof.} We prove {\bf (9)}.
Assume that $f:B_n \to \Z$ is an arithmetic map. Then,
\par
\noindent
$f(1)=1$,
\par
\noindent
$f(9)=f(3 \cdot 3)=f(3) \cdot f(3)=(f(3))^2$,
\par
\noindent
$f(27)=f(3 \cdot 9)=f(3) \cdot f(9)=f(3) \cdot (f(3))^2=(f(3))^3$,
\par
\noindent
$f(25)=f(5 \cdot 5)=f(5) \cdot f(5)=(f(5))^2$,
\par
\noindent
$f(26)=f(25+1)=f(25)+f(1)=(f(5))^2+1$,
\par
\noindent
$f(27)=f(26+1)=f(26)+f(1)=(f(5))^2+1+1=(f(5))^2+2$.
\par
\noindent
Therefore, $(f(3))^3=f(27)=(f(5))^2+2$. The equation $x^3=y^2+2$ has
$(3,\pm 5)$ as its only integer solutions, see [17,~p.~398], [8,~p.~124], [12,~p.~104],
[9,~p.~66], [11,~p.~57].
Thus, $f(3)=3$ and $f(5)=\pm 5$. Hence,
$f(25)=f(5 \cdot 5)=f(5) \cdot f(5)=(\pm 5)^2=25$,
$f(26)=f(25+1)=f(25)+f(1)=25+1=26$.
Since $f(3)=3$, we get by induction $f(3^k)=3^k$ for each $k \in \{1,2,3,...,n\}$.
We have proved~{\bf (9)}. The equation $x^3=y^2+2$ has a rational
solution $(\frac{129}{100},\frac{383}{1000})$, see [2,~p.~173], [13,~p.~2], [9 p.~66],
[11,~p.~57]. \mbox{We define $\phi:B \to \Q$ as}
\begin{eqnarray*}
\left\{\left(1,1\right),\left(5,\frac{383}{1000}\right),\left(25,\left(\frac{383}{1000}\right)^2\right),\left(26,\left(\frac{383}{1000}\right)^2+1\right)\right\}
& \cup & \\
\left\{\left(3,\frac{129}{100}\right) \left(3^2,\left(\frac{129}{100}\right)^2\right),\left(3^3,\left(\frac{129}{100}\right)^3\right),...\right\}
\end{eqnarray*}
The map $\phi$ is arithmetic and $\phi$ moves all $r\in B\setminus \{1\}$, so condition~{\bf (10)} holds true.
\rightline{$\Box$}
\par
We present a simpler counterexample for $r=-1$. Let $G=\{-4,-1,1,3,9,12,16\}$, $\eta: G \to \Q$,
$\eta=\left\{\left(-4,\frac{1}{2}\right),(-1,1),(1,1),\left(3,\frac{1}{2}\right),\left(9,\frac{1}{4}\right),\left(12,\frac{3}{4}\right),\left(16,\frac{1}{4}\right) \right\}$.
The map $\eta$ is arithmetic and $\eta$ moves $-1$.
We prove that $G$ is an arithmetic neighbourhood \mbox{of $-1$} inside~$\Z$.
Assume that $f:G \to \Z$ is an arithmetic map. Since
\par
\noindent
\centerline{$f(-1)=f(-4+3)=f(-4)+f(3)$,}
\par
\noindent
we get $f(-4)=f(-1)-f(3)$. Hence,
\par
\noindent
$f(16)=f((-4) \cdot (-4))=f(-4) \cdot f(-4)=(f(-1)-f(3))^2$,
\par
\noindent
$f(12)=f(-4+16)=f(-4)+f(16)=f(-1)-f(3)+(f(-1)-f(3))^2$.
\par
\noindent
Since $1=f((-1) \cdot (-1))=f(-1) \cdot f(-1)$, we get $f(-1)=\pm 1$.
Assume, on the contrary, that $f(-1)=1$. Thus,
\par
\noindent
\centerline{$1-f(3)+(1-f(3))^2=f(12)=f(3+3 \cdot 3)=f(3)+f(3 \cdot 3)=f(3)+(f(3))^2$}
\par
\noindent
Solving this equation for $f(3)$ we obtain $2 \cdot f(3)=1$, a contradiction.
\vskip 0.2truecm
\par
Let $Y=\{-4,-1,1,3,9,12,14,16,20,180,196\}$, $\kappa:Y \to {\Q}(\sqrt{3})$,
$$
\kappa=\bigl\{\left(-4,\frac{1}{2}\right),(-1,1),(1,1),\left(3,\frac{1}{2}\right),
\left(9,\frac{1}{4}\right),\left(12,\frac{3}{4}\right),
$$
$$
\left(14,\frac{\sqrt{3}}{4}\right),\left(16,\frac{1}{4}\right),\left(20,-\frac{1}{4}\right),
\left(180,-\frac{1}{16}\right),\left(196,\frac{3}{16}\right)\bigr\}
$$
The map $\kappa$ is arithmetic and $\kappa$ moves $-1$.
We prove that $Y$ is an arithmetic neighbourhood \mbox{of $-1$}
inside~$\Q$. Let $f:Y \to \Q$ be an arithmetic map, and
assume, on the contrary, that $f(-1)=1$.
As previously, we conclude that $f(3)=\frac{1}{2}$ and
$f(-4)=f(-1)-f(3)=1-\frac{1}{2}=\frac{1}{2}$.
Hence,
\par
\noindent
$f(9)=f(3 \cdot 3)=f(3) \cdot f(3)=\frac{1}{2} \cdot \frac{1}{2}=\frac{1}{4}$,
\par
\noindent
$f(16)=f((-4) \cdot (-4))=f(-4) \cdot f(-4)=\frac{1}{2} \cdot \frac{1}{2}=\frac{1}{4}$.
\par
\noindent
Since $f(16)=f(-4+20)=f(-4)+f(20)$,
\par
\noindent
we get $f(20)=f(16)-f(-4)=\frac{1}{4}-\frac{1}{2}=-\frac{1}{4}$. Therefore,
\par
\noindent
$f(180)=f(9 \cdot 20)=f(9) \cdot f(20)= \frac{1}{4} \cdot (-\frac{1}{4})=-\frac{1}{16}$,
\par
\noindent
$(f(14))^2=f(14 \cdot 14)=f(180+16)=f(180)+f(16)=-\frac{1}{16}+\frac{1}{4}=\frac{3}{16}$,
a contradiction.
\vskip 0.2truecm
\par
Let $M=\{-4,-1,1,3,5,9,11,42,45,121,126\}$, $\chi: M \to {\Z}[\sqrt{-1}]$,
\par
\noindent
\small
\centerline{$\chi=\{(-4,0),(-1,1),(1,1),(3,1),(5,1),(9,1),(11,\sqrt{-1}),(42,0),(45,1),(121,-1),(126,0)\}$}
\normalsize
\par
\noindent
The map $\chi$ is arithmetic and $\chi$ moves $-1$. We prove
that $M$ is an arithmetic neighbourhood \mbox{of $-1$} inside~$\R$.
Let $f:M \to \R$ be an arithmetic map, and assume, on the contrary,
that $f(-1)=1$. Since $1=f(-4+3)=f(-4)+f(3)$, we get
$f(-4)=1-f(3)$. Hence
\begin{eqnarray*}
1=f(-4+5)=f(-4)+f(5)=1-f(3)+f(-4+3 \cdot 3)
&=&\\
1-f(3)+f(-4)+f(3 \cdot 3)=1-f(3)+1-f(3)+(f(3))^2
\end{eqnarray*}
Solving this equation for $f(3)$ we obtain $f(3)=1$. Therefore,
\par
\noindent
$f(-4)=1-f(3)=1-1=0$,
\par
\noindent
$f(5)=0+f(5)=f(-4)+f(5)=f(-4+5)=1$,
\par
\noindent
$f(9)=f(3 \cdot 3)=f(3) \cdot f(3)=1 \cdot 1=1$,
\par
\noindent
$f(45)=f(5 \cdot 9)=f(5) \cdot f(9)=1 \cdot 1=1$.
\par
\noindent
Since $f(45)=f(3+42)=f(3)+f(42)=1+f(42)$,
we get $f(42)=f(45)-1=1-1=0$.
Thus, $f(126)=f(3 \cdot 42)=f(3) \cdot f(42)=1 \cdot 0=0$.
Since $0=f(5+11 \cdot 11)=f(5)+f(11 \cdot 11)=1+(f(11))^2$,
we get $(f(11))^2=-1$, a contradiction.
\vskip 0.2truecm
\par
\noindent
{\bf Remark 2.} By Theorem 2 infinitely many integers~$r$ do not satisfy condition~{\bf (5)}.
Considering the equation $x^2+2y^2=1025$ one can prove
that the numbers $r=15^2$, $r=20^2$,
$r=2 \cdot 20^2$ do not satisfy condition~{\bf (5)}.
Considering the equation $x^2+y^2=218$ one can prove that $r=7^2 \cdot 13^2$
does not satisfy condition~{\bf (5)}.
Considering the equation $x^2+y^2=1021$ one can prove
that $r=11^2 \cdot 30^2$ does not satisfy condition~{\bf (5)}.
\vskip 0.2truecm
\par
Let $w$ denote the unique real root of the polynomial $x^3-x^2-x-3$.
\vskip 0.2truecm
\par
\noindent
{\bf Theorem 3.} There is an arithmetic map $\psi: \{-4\} \cup B \to {\Q}(w)$
which moves all $r \in \{-4\} \cup B \setminus \{1\}$.
For each $r \in \{-4\} \cup B_n \setminus \{1\}$ we have:
\par
\noindent
{\bf (11)}~~~$\{-4\} \cup B_n$ is an arithmetic neighbourhood of~$r$ inside~$\Q$,
\par
\noindent
{\bf (12)}~~~$\{-4\} \cup B_n$ is not an arithmetic neighbourhood of~$r$ inside~${\Q}(w)$.
\vskip 0.2truecm
\par
\noindent
{\it Proof.} We prove {\bf (11)}.
Assume that $f:\{-4\} \cup B_n \to \Q$ is an arithmetic map. Since
$1=f(1)=f(-4+5)=f(-4)+f(5)$, we get
\begin{equation}
\tag*{{\bf (13)}}
f(-4)=1-f(5)
\end{equation}
\par
\noindent
Hence,
\par
\noindent
$f(5)=f(-4+(3 \cdot 3))=f(-4)+f(3 \cdot 3)=f(-4)+f(3) \cdot f(3)=1-f(5)+(f(3))^2$.
\par
\noindent
Therefore,
\begin{equation}
\tag*{{\bf (14)}}
f(5)=\frac{1+(f(3))^2}{2}
\end{equation}
From equations {\bf (13)} and {\bf (14)}, we obtain
\begin{equation}
\tag*{{\bf (15)}}
f(-4)=1-f(5)=1-\frac{1+(f(3))^2}{2}=\frac{1-(f(3))^2}{2}
\end{equation}
Proceeding exactly as in the proof of Theorem 2, we obtain $(f(3))^3=(f(5))^2+2$.
By this and equation~{\bf (14)}, we get
\begin{equation}
\tag*{{\bf (16)}}
(f(3))^3=\left(\frac{1+(f(3))^2}{2}\right)^2+2
\end{equation}
Equation {\bf (16)} is equivalent to the equation
$$(f(3)-3) \cdot ((f(3))^3-(f(3))^2-f(3)-3)=0$$
The equation $x^3-x^2-x-3=0$ has no rational solutions,
so we must have $f(3)=3$. By induction we get
$f(3^k)=3^k$ for each $k \in \{1,2,3,...,n\}$.
Knowing that $f(3)=3$, from equations {\bf (15)} and {\bf (14)}
we obtain:
$$f(-4)=\frac{1-(f(3))^2}{2}=\frac{1-3^2}{2}=-4$$
$$f(5)=\frac{1+(f(3))^2}{2}=\frac{1+3^2}{2}=5$$
Consequently,
$$f(25)=f(5 \cdot 5)=f(5) \cdot f(5)=5 \cdot 5=25$$
$$f(26)=f(25+1)=f(25)+f(1)=25+1=26$$
The proof of {\bf (11)} is completed.
We define $\psi:\{-4\} \cup B \to {\Q}(w)$ as
\begin{eqnarray*}
\left\{\left(-4,\frac{1-w^2}{2}\right),\left(1,1\right),
\left(5,\frac{1+w^2}{2}\right),\left(25,\left(\frac{1+w^2}{2}\right)^2\right),
\left(26,\left(\frac{1+w^2}{2}\right)^2+1\right)\right\}
& \cup & \\
\left\{\left(3,w\right), \left(3^2,w^2\right),\left(3^3,w^3\right),...\right\}
\end{eqnarray*}
The map $\psi$ is arithmetic and $\psi$ moves all $r \in \{-4\} \cup B \setminus \{1\}$,
so condition~{\bf (12)} holds true.
\newline 
\rightline{$\Box$}
\par
Let $n \in \Z$, $n \geq 1$, $C_n=\left\{1,3,5,13,25,65,169,194,195\right\} \cup \left\{9,9^2,9^3,...,9^n\right\}$,
$C=\bigcup_{n=1}^\infty\limits C_n$.
\vskip 0.2truecm
\par
\noindent
{\bf Theorem 4.} There is an arithmetic map $g: C \to \Q$
which moves all $r \in C \setminus \{1\}$.
\par
\noindent
{\bf (17)}~~~$C_n$ is an arithmetic neighbourhood inside~$\Z$ for $9$, $9^2$, $9^3$, ..., $9^n$,
\par
\noindent
{\bf (18)}~~~$C_n$ is not an arithmetic neighbourhood inside~$\Q$ for $9$, $9^2$, $9^3$, ..., $9^n$.
\vskip 0.2truecm
\par
\noindent
{\it Proof.} We prove {\bf (17)}.
Assume that $f:C_n \to \Z$ is an arithmetic map. Then,
\vskip 0.2truecm
\par
\noindent
$(f(5))^2+(f(13))^2+1=f(5^2)+f(13^2)+f(1)=f(5^2+13^2)+f(1)=f(5^2+13^2+1)=$
\vskip 0.2truecm
\par
\noindent
\centerline{$f((5 \cdot 13) \cdot 3)=f(5 \cdot 13) \cdot f(3)=f(5) \cdot f(13) \cdot f(3)$.}
\vskip 0.2truecm
\par
\noindent
If integers $x$, $y$, $z$ satisfy $x^2+y^2+1=xyz$ then $z=\pm 3$,
see [8,~p.~299], [9,~pp.~58--59], [10,~p.~31], [11,~pp.~51--52], \cite{Barnes}, cf. Theorem~4 in [7,~p.~218].
Thus, $f(3)=\pm 3$. Hence $f(9)=f(3 \cdot 3)=f(3) \cdot f(3)=(\pm 3)^2=9$,
and by induction we obtain $f(9^k)=9^k$ for each $k \in \{1,2,3,...,n\}$.
The proof of~{\bf (17)} is completed.
We define $g:C \to \Q$ as
\begin{eqnarray*}
\left\{(1,1),\left(3,\frac{9}{4}\right),(5,2),(13,2),(25,4),(65,4),(169,4),(194,8), (195,9)\right\}
& \cup & \\
\left\{\left(9,\frac{81}{16}\right),
\left(9^2,\left(\frac{81}{16}\right)^2\right),
\left(9^3,\left(\frac{81}{16}\right)^3\right),...,
\left(9^n,\left(\frac{81}{16}\right)^n\right)
\right\}
\end{eqnarray*}
The map $g$ is arithmetic and $g$ moves all $r \in C \setminus \{1\}$,
so condition~{\bf (18)} holds true.
\newline
\rightline{$\Box$}
\par
We know (see Theorem 2 or Theorem 4) that infinitely many integers~$r$
do not satisfy condition~{\bf (5)}. Now, we sketch a more elementary
(but longer) proof of this fact. Let $n \in \Z$, $n \ge 3$,
\begin{eqnarray*}
H_n=\{1,2,4,16,60,64,3600,3604,3620,3622,3623,7^3 \cdot 13^2,7^3 \cdot 13^2+1\}
& \cup & \\
\{13,13^2,7,7^2,7^3, ..., 7^n\}
\end{eqnarray*}
$H_n$ is an arithmetic neighbourhood inside $\Z$
for each $r \in H_n \setminus \{13\}$. $H_n$ is not an
arithmetic neighbourhood inside $\Q$ for
$13$, $13^2$, $7$, $7^2$, $7^3$, ..., $7^n$.
The proofs follow from the following observations:
\par
\noindent
\centerline{$7^3 \cdot 13^2+1=16 \cdot 3623$}
\par
\noindent
\centerline{$\forall x,y \in \Z ~(x^3 \cdot y^2=7^3 \cdot 13^2 \Rightarrow (x=7 \wedge y=\pm 13))$}
\par
\noindent
\centerline{$(\frac{7}{4})^3 \cdot (8 \cdot 13)^2=7^3 \cdot 13^2$}
\newpage
\par
\noindent
{\bf Theorem 5.} Let
\par
\noindent
\centerline{$D=\left\{-36,\frac{1}{2},1,2,\frac{5}{2},5,12,25,50,100,12^2,200,400,425,430,432,36^2,12^3\right\}$.}
\par
\noindent
{\bf (19)}~~~$D$ is an arithmetic neighbourhood inside~$\Q$ for $12$, $12^2$, $36^2$, $12^3$,
\par
\noindent
{\bf (20)}~~~$D$ is not an arithmetic neighbourhood inside~${\Q}(\sqrt{5})$ for $12$, $12^2$, $36^2$, $12^3$.
\vskip 0.2truecm
\par
\noindent
{\it Proof.} We prove {\bf (19)}. Assume that $f:D \to \Q$ is an arithmetic map.
Then, $f(1)=1$ and $f(2)=f(1+1)=f(1)+f(1)=1+1=2$. Since
$1=f(1)=f(\frac{1}{2}+\frac{1}{2})=f(\frac{1}{2})+f(\frac{1}{2})$,
we get $f(\frac{1}{2})=\frac{1}{2}$. Knowing $f(\frac{1}{2})$ and $f(2)$, we calculate
\par
\noindent
$f(\frac{5}{2})=f(2+\frac{1}{2})=f(2)+f(\frac{1}{2})=2+\frac{1}{2}=\frac{5}{2}$,
\par
\noindent
$f(5)=f(\frac{5}{2}+\frac{5}{2})=f(\frac{5}{2})+f(\frac{5}{2})=\frac{5}{2}+\frac{5}{2}=5$,
\par
\noindent
$f(25)=f(5 \cdot 5)=f(5) \cdot f(5)=5 \cdot 5=25$,
\par
\noindent
$f(50)=f(25+25)=f(25)+f(25)=25+25=50$,
\par
\noindent
$f(100)=f(50+50)=f(50)+f(50)=50+50=100$,
\par
\noindent
$f(200)=f(100+100)=f(100)+f(100)=100+100=200$,
\par
\noindent
$f(400)=f(200+200)=f(200)+f(200)=200+200=400$,
\par
\noindent
$f(425)=f(400+25)=f(400)+f(25)=400+25=425$,
\par
\noindent
$f(430)=f(425+5)=f(425)+f(5)=425+5=430$,
\par
\noindent
$f(432)=f(430+2)=f(430)+f(2)=430+2=432$.
\vskip 0.2truecm
\par
\noindent
Therefore, $(f(12))^3=(f(12) \cdot f(12)) \cdot f(12)=f(12 \cdot 12) \cdot f(12)=
f((12 \cdot 12) \cdot 12)=f((-36)^2+432)=f((-36)^2)+f(432)=(f(-36))^2+432$.
The equation \mbox{$x^3=y^2+432$} has $(12,\pm 36)$ as its only rational solutions,
see \cite{Fueter}, [12,~p.~107], [2,~p.~174], [4,~p.~296], [8, p.~247], [14,~p.~54].
Thus, $f(12)=12$ and $f(-36)=\pm 36$.
Hence, $f(12^2)=f(12) \cdot f(12)=12^2$,
$f(12^3)=f(12 \cdot 12^2)=f(12) \cdot f(12^2)=12 \cdot 12^2=12^3$,
$f(36^2)=f((-36) \cdot (-36))=f(-36) \cdot f(-36)=(\pm 36)^2=36^2$.
The proof of~{\bf (19)} is completed.
We find that \mbox{$8^3=(4 \cdot\sqrt{5})^2+432$} and we define
$h:D \to {\Q}(\sqrt{5})$ as
\begin{eqnarray*}
\left\{(-36,4 \cdot \sqrt{5}),(12,8),(12^2,8^2),(36^2,80), (12^3,8^3)\right\} & \cup & \\
{\rm id} \left( \left\{ \frac{1}{2},1,2,\frac{5}{2},5,25,50,100,200,400,425,430,432 \right\} \right)
\end{eqnarray*}
We summarize the check that $h$ is arithmetic. Obviously, $h(1)=1$.
To check the condition
\par
\noindent
\centerline{$\forall x,y,z \in D ~(x+y=z \Rightarrow h(x)+h(y)=h(z))$}
\par
\noindent
it is enough to consider all the triples $(x,y,z) \in D \times D \times D$
for which $x+y=z$, $x \leq y$, and $h$ is not the identity on $\left\{x,y,z\right\}$.
There is only one such triple: $(432,36^2,12^3)$.
\newpage
\par
\noindent
To check the condition
\par
\noindent
\centerline{$\forall x,y,z \in D ~(x \cdot y=z \Rightarrow h(x) \cdot h(y)=h(z))$}
\par
\noindent
it is enough to consider all the triples $(x,y,z) \in D \times D \times D$
for which $x \cdot y=z$, $x \leq y$, $x \neq 1$, $y \neq 1$, and $h$
is not the identity on $\left\{x,y,z\right\}$. These triples are as follows:
\par
\noindent
\centerline{$(-36,-36,36^2)$, $(12,12,12^2)$, $(12,12^2,12^3)$}
\par
\noindent
The sentence~{\bf (20)} is true because $h$ is arithmetic and $h$ moves
$12$, $12^2$, $36^2$, $12^3$.
\newline
\rightline{$\Box$}
\par
\noindent
{\bf Corollary.} Let us define by induction the finite sets
$D_n \subseteq \Q$ ($n=0,1,2,...$). Let $D_0=D$, $d_n$~denote
the greatest number in $D_n$, $D_{n+1}=D_n \cup \left\{d_n^2\right\}$.
For each $n \in \left\{0,1,2,...\right\}$ we have:
\par
\noindent
$D_n$ is an arithmetic neighbourhood of $d_n$ inside~$\Q$,
\par
\noindent
$D_n$ is not an arithmetic neighbourhood of $d_n$ inside~${\Q}(\sqrt{5})$.
\vskip 0.2truecm
\par
Let $u=\frac{-1 \pm \sqrt{33}}{8}$, $n \in \Z$, $n \geq 3$,
$E_n=\left\{\frac{1}{2},1,\frac{3}{2},\frac{9}{4}\right\} \cup \left\{9,-2,(-2)^2,(-2)^3,...,(-2)^n\right\}$,
$E=\bigcup_{n=3}^\infty\limits E_n$.
\vskip 0.2truecm
\par
\noindent
{\bf Theorem 6.} There is an arithmetic map $\sigma: E \to {\Q}(\sqrt{33})$
which moves $9$ and all the numbers $(-2)^k$, where $k \in \{1,2,3,...\}$.
For each \mbox{$r \in E_n \setminus \{\frac{1}{2},1,\frac{3}{2},\frac{9}{4}\}$}
we have:
\par
\noindent
{\bf (21)}~~~$E_n$ is an arithmetic neighbourhood of $r$ inside~$\Q$,
\par
\noindent
{\bf (22)}~~~$E_n$ is not an arithmetic neighbourhood of $r$ inside~${\Q}(\sqrt{33})$.
\vskip 0.2truecm
\par
\noindent
{\it Proof.} We prove {\bf (21)}.
Assume that $f:E_n \to \Q$ is an arithmetic map.
\vskip 0.2truecm
\par
\noindent
Since $1=f(\frac{1}{2}+\frac{1}{2})=f(\frac{1}{2})+f(\frac{1}{2})$,
we get $f(\frac{1}{2})=\frac{1}{2}$. Hence,
\par
\noindent
$f(\frac{3}{2})=f(1+\frac{1}{2})=f(1)+f(\frac{1}{2})=1+\frac{1}{2}=\frac{3}{2}$. Thus,
\par
\noindent
$f(\frac{9}{4})=f(\frac{3}{2} \cdot \frac{3}{2})
=f(\frac{3}{2}) \cdot f(\frac{3}{2})=\frac{3}{2} \cdot \frac{3}{2}=\frac{9}{4}$.
Therefore,
\par
\noindent
$f(9)=f(\frac{9}{4} \cdot 4)=f(\frac{9}{4}) \cdot f((-2) \cdot (-2))=
\frac{9}{4} \cdot (f(-2))^2$. It implies that
\begin{eqnarray*}
1=f(-2 \cdot 4+9)=f(-2 \cdot 4)+f(9)=
f(-2) \cdot f(4)+\frac{9}{4} \cdot (f(-2))^2
&=&\\
f(-2) \cdot f(-2 \cdot (-2))+\frac{9}{4} \cdot (f(-2))^2=
(f(-2))^3+\frac{9}{4} \cdot (f(-2))^2
\end{eqnarray*}
Solving this equation for $f(-2)$ we obtain $f(-2)=-2$,
the only rational root.
Another roots are $\frac{-1-\sqrt{33}}{8}$ and $\frac{-1+\sqrt{33}}{8}$.
Knowing $f(-2)$, we calculate
\par
\noindent
\centerline{$f(9)=\frac{9}{4} \cdot (f(-2))^2=\frac{9}{4} \cdot (-2)^2=9$}
\par
\noindent
Applying induction, we obtain $f((-2)^k)=(-2)^k$ for each $k \in \{1,2,3,...,n\}$.
We have proved~{\bf (21)}. We define $\sigma:E \to {\Q}(\sqrt{33})$ as
$$
{\rm id}\left(\left\{\frac{1}{2},1,\frac{3}{2},\frac{9}{4}\right\}\right) \cup
\left\{\left(9,\frac{9}{4} \cdot u^2 \right),\left(-2,u \right),\left((-2)^2,u^2\right),\left((-2)^3,u^3\right),~...\right\}
$$
The map $\sigma$ is arithmetic and $\sigma$ moves all
\mbox{$r \in E_n \setminus \{\frac{1}{2},1,\frac{3}{2},\frac{9}{4}\}$},
so condition~{\bf (22)} holds true.
\newline
\rightline{$\Box$}
\vskip 0.2truecm
\par
Theorem 7 which follows is more general than the previous ones.
Let $n$ be an integer, and assume that $n \geq 3$ and
$n \not\in \{2^2,2^3,2^4,...\}$.
We find the smallest integer $\rho(n)$ such that $n^3 \leq 2^{\rho(n)}$.
From the definition of $\rho(n)$ we obtain $2^{\rho(n)-1}<n^3$. It gives
$$
2^{\rho(n)}=2 \cdot 2^{\rho(n)-1}<2 \cdot n^3<n \cdot n^3=n^4
$$
Since $n^3 \leq 2^{\rho(n)}<n^4$, $2^{\rho(n)}$ has four digits in the number system with base~$n$. Let
$$ 
2^{\rho(n)}=m_3 \cdot n^3+m_2 \cdot n^2+m_1 \cdot n+m_0
$$
where $m_3 \in \{1,2,...,n-1\}$ and $m_2,m_1,m_0 \in \{0,1,2,...,n-1\}$. Let
\begin{eqnarray*}
{\cal J}(n)=\left\{-1,~0,~1,~-\frac{1}{2},~-\frac{1}{2^2},~-\frac{1}{2^3},~...,~-\frac{1}{2^{\rho(n)}},~n,~n^2 \right\} \cup \\
\left\{k \cdot n^3:~~ k \in \{1,2,...,m_3\}\right\}\cup\\
\{m_3 \cdot n^3+k \cdot n^2:~~ k \in \{1,2,...,m_2\}\}\cup\\
\{m_3 \cdot n^3+m_2 \cdot n^2+k \cdot n:~~ k \in \{1,2,...,m_1\}\}\cup\\
\{m_3 \cdot n^3+m_2 \cdot n^2+m_1 \cdot n+k:~~ k \in \{1,2,...,m_0\}\}~~~
\end{eqnarray*}
Of course, $\{2^{\rho(n)},~2^{\rho(n)}-1,~2^{\rho(n)}-2,~...,~2^{\rho(n)}-m_0\} \subseteq {\cal J}(n)$.
\vskip 0.2truecm
\par
\noindent
{\bf Theorem 7}. ${\cal J}(n)$ is an arithmetic neighbourhood of~$n$
inside~$\R$, and so too inside~$\Q$. ${\cal J}(n)$ is not an arithmetic
neighbourhood of~$n$ inside~$\C$.
\vskip 0.2truecm
\par
\noindent
{\it Proof.} Assume that $f: {\cal J}(n) \to \R$ is an arithmetic
map. Since $0+0=0$, $f(0)=0$.
Since $f(0)=0$ and $-1+1=0$, $f(-1)=-1$. Since
$$
-1=f\left(-1\right)=f\left(\left(-\frac{1}{2}\right)+\left(-\frac{1}{2}\right)\right)=f\left(-\frac{1}{2}\right)+f\left(-\frac{1}{2}\right)
$$
we get $f\left(-\frac{1}{2}\right)=-\frac{1}{2}$. Hence,
from $-\frac{1}{2}=\left(-\frac{1}{4}\right)+\left(-\frac{1}{4}\right)$
we obtain $f\left(-\frac{1}{4}\right)=-\frac{1}{4}$. Applying induction we obtain
$f\left(-\frac{1}{2^{\rho(n)}}\right)=-\frac{1}{2^{\rho(n)}}$. We have
$$
-1=f\left(-1\right)=f\left(2^{\rho(n)} \cdot ~\left(-\frac{1}{2^{\rho(n)}}\right)\right)=
f\left(2^{\rho(n)}\right) \cdot f\left(-\frac{1}{2^{\rho(n)}}\right)=
f(2^{\rho(n)}) \cdot ~\left(-\frac{1}{2^{\rho(n)}}\right)
$$
Hence $f(2^{\rho(n)})=2^{\rho(n)}$. Applying induction we get
$$
2^{\rho(n)}=f(2^{\rho(n)})=f(m_3 \cdot n^3+m_2 \cdot n^2+m_1 \cdot n+m_0)=m_3 \cdot (f(n))^3+m_2 \cdot (f(n))^2+m_1 \cdot f(n)+m_0
$$
We want to prove that $f(n)=n$. It suffices to show that the function
$$
\R \ni x \stackrel{\zeta}{\longrightarrow} m_3 \cdot x^3+m_2 \cdot x^2+m_1 \cdot x+m_0 \in \R
$$
takes the value $2^{\rho(n)}$ only for $x=n$.
Since $\zeta$ is strictly increasing in the interval $[0,\infty)$,
for each $x \in [0,n)$ we have $\zeta(x)<\zeta(n)=2^{\rho(n)}$, and for each $x \in (n,\infty)$
we have $\zeta(x)>\zeta(n)=2^{\rho(n)}$. We show that $\zeta$ does not
reach the value $2^{\rho(n)}$ for $x \in (-\infty,0]$.
For each $x \in (-\infty,0]$ we have
\begin{equation}
\tag*{{\bf (23)}}
\zeta(x)=m_3 \cdot x^3+m_2 \cdot x^2+m_1 \cdot x+m_0 \leq x^3+(n-1) \cdot x^2+n-1
\end{equation}
By {\bf (23)}, if $x \in (-\infty,-n+1]$ then
$$
\zeta(x) \leq x^3+(n-1) \cdot x^2+n-1=(x+n-1) \cdot x^2+n-1 \leq n-1<n^3 \leq 2^{\rho(n)}
$$
Thus, $\zeta(x) \neq 2^{\rho(n)}$. By {\bf (23)}, if $x \in [-n+1,0]$ then
$$
\zeta(x) \leq x^3+(n-1) \cdot x^2+n-1 \leq (n-1) \cdot x^2+n-1 \leq (n-1)^3+n-1<n^3 \leq 2^{\rho(n)}
$$
Thus, $\zeta(x) \neq 2^{\rho(n)}$.
We have proved that $f(n)=n$. It proves that ${\cal J}(n)$ is an arithmetic
neighbourhood of~$n$ inside~$\R$. We prove that ${\cal J}(n)$
is not an arithmetic neighbourhood of $n$ inside $\C$.
The number~$n$ is a single root of the polynomial
$$
m_3 \cdot x^3+m_2 \cdot x^2+m_1 \cdot x+m_0-2^{\rho(n)}
$$
because the derivative of this polynomial takes the non-zero value
$$
3 \cdot m_3 \cdot n^2+2 \cdot m_2 \cdot n+m_1 \geq 3 \cdot n^2 \geq 27
$$
at $x=n$. Hence the polynomial
$$
m_3 \cdot x^3+m_2 \cdot x^2+m_1 \cdot x+m_0-2^{\rho(n)}
$$
has two conjugated roots $z_1,z_2 \in \C \setminus \R$.
Let $z=z_1$ or $z=z_2$. We define $\theta: {\cal J}(n) \to \C$ as
\begin{eqnarray*}
{\rm id}\left(\left\{-1,~0,~1,~-\frac{1}{2},~-\frac{1}{2^2},~-\frac{1}{2^3},~...,~-\frac{1}{2^{\rho(n)}}\right\}\right) \cup \{(n,z),~(n^2,z^2)\} \cup \\
\{(k \cdot n^3,~k \cdot z^3):~~k \in \{1,2,...,m_3\}\} \cup \\
\{(m_3 \cdot n^3+k \cdot n^2,~m_3 \cdot z^3+k \cdot z^2):~~k \in \{1,2,...,m_2\}\} \cup \\
\{(m_3 \cdot n^3+m_2 \cdot n^2+k \cdot n,~m_3 \cdot z^3+m_2 \cdot z^2+k \cdot z):~~k \in \{1,2,...,m_1\}\} \cup \\
\{(m_3 \cdot n^3+m_2 \cdot n^2+m_1 \cdot n+k,~m_3 \cdot z^3+m_2 \cdot z^2+m_1 \cdot z+k):~~k \in \{1,2,...,m_0\}\}
\end{eqnarray*}
\vskip 0.2truecm
\par
\noindent
Of course, $\theta(x)=x$ for each $x \in \{2^{\rho(n)},~2^{\rho(n)}-1,~2^{\rho(n)}-2,~...,~2^{\rho(n)}-m_0\}$.
Since $\theta(n)=z \neq n$, $\theta$ moves $n$.
We summarize the check that $\theta$ is arithmetic. Obviously, $\theta(1)=1$.
To check the condition
\par
\noindent
\centerline{$\forall x,y,z \in {\cal J}(n) ~(x+y=z \Rightarrow \theta(x)+\theta(y)=\theta(z))$}
\par
\noindent
it is enough to consider all the triples $(x,y,z) \in {\cal J}(n) \times {\cal J}(n) \times {\cal J}(n)$
for which $x+y=z$, $x \leq y$, $x \neq 0$, $y \neq 0$, and $\theta$ is not the identity on $\left\{x,y,z\right\}$.
These triples are as follows:
\par
\noindent
$(k \cdot n^3,~l \cdot n^3,~(k+l) \cdot n^3)$, where $k,l,k+l \in \{1,2,...,m_3\}$ and $k \leq l$,
\par
\noindent
$(n^2,~m_3 \cdot n^3+k \cdot n^2,~m_3 \cdot n^3+(k+1) \cdot n^2)$, where $k,k+1 \in \{0,1,2,...,m_2\}$,
\par
\noindent
$(n,~m_3 \cdot n^3+m_2 \cdot n^2+k \cdot n,~m_3 \cdot n^3+m_2 \cdot n^2+(k+1) \cdot n)$, where $k,k+1 \in \{0,1,2,...,m_1\}$.
\par
\noindent
To check the condition
\par
\noindent
\centerline{$\forall x,y,z \in {\cal J}(n) ~(x \cdot y=z \Rightarrow \theta(x) \cdot \theta(y)=\theta(z))$}
\par
\noindent
it is enough to consider all the triples $(x,y,z) \in {\cal J}(n) \times {\cal J}(n) \times {\cal J}(n)$
for which $x \cdot y=z$, $x \leq y$, $x \neq 1$, $y \neq 1$, $x \neq 0$, $y \neq 0$, and $\theta$
is not the identity on $\left\{x,y,z\right\}$. These triples are as follows:~
$(n,n,n^2)$, $(n,n^2,n^3)$.
\newline
\rightline{$\Box$}

Apoloniusz Tyszka\\
Technical Faculty\\
Hugo Ko\l{}\l{}\k{a}taj University\\
Balicka 116B, 30-149 Krak\'ow, Poland\\
E-mail address: {\it rttyszka@cyf-kr.edu.pl}

\begin{thebibliography}{16}
\bibitem{Barnes}
{\sc E.~S.~Barnes}, {\it On the Diophantine equation $x^2+y^2+c=xyz$},
J. London Math. Soc. 28 (1953), 242--244.
\bibitem{Bundschuh}
{\sc P.~Bundschuh}, {\it Einf\"uhrung in die {Z}ahlentheorie},
5., \"uberarbeitete und aktualisierte Aufl., Springer, Berlin, 2002.
\bibitem{Fueter}
{\sc R.~Fueter}, {\it Ueber kubische diophantische Gleichungen},
Comment. Math. Helv. 2~(1930), no.~1, 69--89.
\bibitem{Ireland-Rosen}
{\sc K.~Ireland and M.~Rosen},
{\it A classical introduction to modern number theory}, 2nd ed.,
Springer, New-York, 1990.
\bibitem{Lettl}
{\sc G.~Lettl}, {\it Finitely arithmetically fixed elements of a field},
Arch. Math. (Basel) 87 (2006), no.~6, 530--538.
\bibitem{Marcja}
{\sc A.~Marcja and C.~Toffalori},
{\it A guide to classical and modern model theory},
Kluwer Academic Publishers, Dordrecht, 2003.
\bibitem{Mills}
{\sc W.~H.~Mills}, {\it A system of quadratic Diophantine equations},
Pacific J. Math. 3~(1953), no.~1, 209--220.
\bibitem{Mordell}
{\sc L.~J. Mordell}, {\it Diophantine equations}, Academic Press, London, 1969.
\bibitem{Sierpinski1956}
{\sc W.~Sierpi{\'n}ski}, {\it On solution of equations in integers} (Polish),
PWN (Polish Scientific Publishers), Warsaw, 1956.
\bibitem{Sierpinski1959}
{\sc W.~Sierpi{\'n}ski}, {\it Theory of numbers, Part 2} (Polish),
PWN (Polish Scientific Publishers), Warsaw, 1959.
\bibitem{Sierpinski1961}
{\sc W.~Sierpi{\'n}ski}, {\it On solution of equations in integers} (Russian, translated from Polish),
Gos. Izd. Fiz. Mat. Lit., Moscow, 1961.
\bibitem{Sierpinski1987}
{\sc W.~Sierpi{\'n}ski}, {\it Elementary theory of numbers}, 2nd~ed. (ed.~A.~Schinzel),
PWN (Polish Scientific Publishers) and North-Holland, Warsaw-Amsterdam, 1987.
\bibitem{Silverman}
{\sc J.~H.~Silverman and J.~Tate}, {\it Rational points on elliptic curves},
Springer, New-York, 1992.
\bibitem{Toth}
{\sc G.~Toth}, {\it Glimpses of algebra and geometry}, 2nd ed.,
Springer, New-York, 2002.
\bibitem{Tyszka1}
{\sc A.~Tyszka}, {\it A discrete form of the theorem that each field
endomorphism of}~~$\R$~(${\Q}_p$) {\it is the identity},
Aequationes Math. 71~(2006), no.~1--2, 100--108.
\bibitem{Tyszka2}
{\sc A.~Tyszka}, {\it On $\emptyset$-definable elements in a field},
Collect. Math. 58 (2007), no.~1, 73--84.
\bibitem{Uspensky}
{\sc J.~V.~Uspensky and M.~A.~Heaslet}, {\it Elementary number theory},
McGraw-Hill, New-York, 1939.
\end{thebibliography}
\end{document}